\documentclass[11pt]{article}
\usepackage{amssymb}
\newtheorem{lem}{Lemma}[section]

\newtheorem{pro}{Proposition}[section]
\newtheorem{thm}{Theorem}[section]

\makeatletter\@addtoreset{equation}{section}
\renewcommand\theequation{\thesection.\@arabic\c@equation}

\newenvironment{Literature}[1]
{\begin{thebibliography}{#1}}{\end{thebibliography}}

\begin{document}

\begin{center}
{\LARGE \bf    Some  extensions of the  Einstein-Dirac Equation }

\bigskip  \bigskip
{\large  Eui Chul Kim}
\end{center}

\bigskip   \noindent
Department of Mathematics, College of Education, \\   Andong
National University,
Andong 760-749, South  Korea     \\
e-mail: eckim@andong.ac.kr

\bigskip \bigskip \noindent
{\bf Abstract:}  We considered an extension of the standard
functional for the Einstein-Dirac equation  where the Dirac
operator is replaced by the square of the Dirac operator and a
real parameter controlling the length of spinors is introduced.
For one distinguished value of the parameter, the resulting
Euler-Lagrange equations provide a new type of Einstein-Dirac
coupling. We establish a special method for constructing global
smooth solutions of a newly derived Einstein-Dirac system called
the {\it CL-Einstein-Dirac equation of type II} (see Definition
3.1).


\bigskip  \bigskip \bigskip
\noindent
\section{Introduction}

\indent Let $(Q^{ n , r} , \eta)$ be an n-dimensional (connected
smooth) pseudo-Riemannian manifold, where
 the index $r$  is the number of negative eigenvalues of the metric $\eta$.
Assume that $(Q^{n , r} , \eta)$  is  space- and time-oriented and
has a fixed spin structure [1].  For simplicity, we will often
write $Q$ to mean $Q^{n, r}$.   Let $\Sigma (Q) = \Sigma
(Q)_{\eta}$ denote the spinor bundle of $(Q^{n, r} , \eta)$
equipped with the ${\rm Spin}^+ (n , r)$-equivariant nondegenerate
complex product $\langle  \cdot ,  \cdot  \rangle = \langle  \cdot
,  \cdot \rangle_{\eta}$, and  let  $( \cdot  ,  \cdot  )  = {\rm
Re} \langle  \cdot ,  \cdot  \rangle $  denote the real part of
$\langle \cdot ,  \cdot  \rangle$. Let ${\rm Ric} = {\rm
Ric}_{\eta}$ and $S = S_{\eta}$ be the Ricci tensor and the scalar
curvature of $(Q^{n, r}, \eta)$, respectively. Let $D = D_{\eta}$
be the Dirac operator acting  on sections  $\psi   \in  \Gamma(
\Sigma (Q) )$ of the spinor bundle $\Sigma (Q)$.  Then the
standard functional for the Einstein-Dirac equation is given by
\begin{equation}
W_1 ( \eta, \psi ) =  \int   \Big\{ a  S_{\eta}  + b  +  \epsilon
\nu_1  (\psi, \psi )  -  \epsilon  (  ( \sqrt{-1} )^r D_{\eta}
\psi,  \psi )    \Big\}   \mu_{\eta},
\end{equation}
where $a, b, \epsilon,  \nu_1   \in {\mathbb R}$, $\epsilon  \neq
0$, are real numbers and $\mu_{\eta}$ is the volume form of $(
Q^{n, r}, \eta)$. The Euler-Lagrange equations (called the
Einstein-Dirac equation) are the Dirac equation
\begin{equation}
 ( \sqrt{-1} )^r  D  \psi    \  =  \    \nu_1 \,   \psi
 \end{equation}
and the Einstein equation
\begin{equation}
 a   \Big\{  {\rm Ric} -   \frac{S}{2}  \,     \eta  \Big\}  - \frac{b}{2}  \,  \eta =  \frac{\epsilon}{4}  \,
 T_1
\end{equation}
coupled via a symmetric tensor field $T_1$,
\begin{equation}
T_1 (X , Y) =  \Big( {(\sqrt{-1})}^r  \{ X \cdot \nabla_Y \psi  +
Y \cdot \nabla_X \psi \}  , \, \psi  \Big) ,
\end{equation}
where $X, Y$  are vector fields on $Q^{n,r}$ and the dot $" \cdot
"$ indicates the Clifford multiplication.  Observe that the system
(1.2)-(1.4) contains four differential operators, namely, the spin
connection $\nabla$, the Dirac operator $D$, the Ricci tensor
${\rm  Ric}$ and the scalar curvature $S$.  The spin connection
and the Dirac operator act on spinor fields and are operators of
first-order,
 while the Ricci tensor and the scalar curvature are second-order operators acting on metrics.
Therefore, it  is natural  to ask whether  one can derive such
Euler-Lagrange equations from the functional
\begin{equation}
W_2 ( \eta, \psi ) =  \int   \Big\{ a  S_{\eta}  + b  +  \epsilon
\nu_2  (\psi, \psi )  -  \epsilon  (  ( D_{\eta}  \circ D_{\eta} )
( \psi) ,  \psi )    \Big\}   \mu_{\eta},   \qquad   \nu_2  \in
{\mathbb R},
\end{equation}
that   generalize the system (1.2)-(1.4) and all the involved
operators acting on spinor fields are   of second-order. In
Section 2 we will show that the answer of the question is positive
and (1.5) yields in fact the following system (see Theorem 2.1):
\begin{equation}
  D^2 \psi =   \nu_2  \psi,         \qquad
 a \Big\{ {\rm Ric}  -  \frac{S}{2}  \,  \eta \Big\} -  \frac{b}{2}  \,  \eta  =   \frac{\epsilon}{4}  T_2   ,
\end{equation}   where $T_2$ is a symmetric tensor field  defined by
\begin{eqnarray}
T_2  (X, Y)   &   =   &    \Big( X \cdot  \nabla_Y (D \psi) +  Y \cdot  \nabla_X (D  \psi) ,   \    \psi  \Big)          \nonumber        \\
&          &       \nonumber        \\
 &            &    +   (-1)^r  \Big( X \cdot  \nabla_Y \psi +  Y \cdot  \nabla_X  \psi ,   \   D  \psi
 \Big).
\end{eqnarray}
In this paper the system (1.2)-(1.4) is called the {\it classical
Einstein-Dirac equation of type I} [5, 6, 7]  and the system
(1.6)-(1.7) the {\it classical Einstein-Dirac equation of type
II}.

\par
Let us turn to another situation where a real parameter
controlling the length of spinors is introduced.  Let $\varphi=
\varphi_{\eta}$ be a spinor field on $(Q^{n,r}, \eta)$ such that
either $(\varphi, \varphi) > 0$ at all points  or $(\varphi,
\varphi) < 0$ at all points.  Fix a shorthand notation
\[     \varphi^k  :=  ( \sigma  \varphi, \varphi )^k   \varphi,   \qquad  \varphi^0 := \varphi,  \]
where $k \in  {\mathbb  R}$ is a real number and  $\sigma =
\sigma_{\varphi} \in {\mathbb R}$ is a constant defined by
\[    \sigma = 1 \   \mbox{if} \  (\varphi, \varphi) >0  \quad  \mbox{and}  \quad   \sigma = -1  \  \mbox{if} \    (\varphi, \varphi) < 0.     \]
Combining the functional (1.1) with (1.5), we extend the spinorial
part as
\begin{equation}
  W ( \eta ,  \varphi )
    =      \int   \Big\{ a  S_{\eta}   + b +  \epsilon    \nu     ( \sigma  \, \varphi^k ,  \varphi^k  )  -   \epsilon   ( \sigma \, P_{\eta} (\varphi^k) , \, \varphi^k )
\Big\}   \mu_{\eta}    ,    \qquad    \nu  \in  {\mathbb R},
\end{equation}
where $P_{\eta} = ( \sqrt{-1} )^r  D_{\eta}$ or $P_{\eta} =
D_{\eta}  \circ D_{\eta}$, and look at the Euler-Lagrange
equations derived from (1.8). We will show in Section 3 (see
Theorem 3.1) that, when $k \not= - \frac{1}{2}$, the
Euler-Lagrange equations of (1.8) are actually equivalent to the
system (1.2)-(1.4) or to the system (1.6)-(1.7) depending on a
choice of $P_{\eta}$.  However, in the distinguished case $k = -
\frac{1}{2}$ in which  the length $\vert \varphi^k \vert =  \pm 1$
becomes constant, we are led to a new Einstein-Dirac system, i.e.,
  \begin{equation}          P_{\eta}  \psi= f  \psi,    \qquad     a  \Big\{ {\rm Ric}  -  \frac{S}{2}    \eta  \Big\}  -    \frac{c}{2} \eta  =   \frac{\epsilon}{4}  T   -  \frac{\epsilon}{2} \, f \, \eta
  ,  \qquad  a , c , \epsilon  \in  {\mathbb R} ,
\end{equation}
where  $\psi$ is of constant length $\vert \psi \vert = \pm 1$ and
$f : Q^{n,r} \longrightarrow {\mathbb R}$ is a real-valued
function and $T$ is a symmetric tensor field defined by
\begin{equation}
 T (X , Y) =  \Big( \sigma {(\sqrt{-1})}^r  \{ X \cdot \nabla_Y \psi  +  Y \cdot \nabla_X \psi \}  , \, \psi  \Big)
\end{equation}
if $P_{\eta} = ( \sqrt{-1} )^r  D_{\eta}$ and  by
\begin{eqnarray}
T (X , Y) &  =  & \sigma \Big(  X \cdot \nabla_Y (D \psi)  +  Y \cdot \nabla_X (D \psi)  , \, \psi  \Big)     \nonumber      \\
&           &        \nonumber       \\
&           &   +  \sigma  (-1)^r  \Big( X \cdot \nabla_Y  \psi +
Y \cdot  \nabla_X \psi ,   \,  D \psi \Big)
\end{eqnarray}
if $P_{\eta} = D_{\eta}  \circ D_{\eta}$, respectively.  The
system (1.9)-(1.11) will be called the {\it CL-Einstein-Dirac
equation of type I} if $P_{\eta} = ( \sqrt{-1} )^r  D_{\eta}$ and
the {\it CL-Einstein-Dirac equation of type II} if $P_{\eta} =
D_{\eta} \circ D_{\eta}$, respectively ("CL" means the "constant
length" of spinors).  A non-trivial spinor field $\psi$ on $(
Q^{n,r}, \eta)$  is called  a {\it CL-Einstein spinor} of type I
(resp. type II) if it satisfies the  CL-Einstein-Dirac equation of
type I (resp. type II).  It will be pointed out (see Remark 3.1)
why one can not weaken the "constant length" condition for
CL-Einstein spinors.

\par
Sections 4 and 5 of the paper are devoted to establishing a
special method for constructing global (smooth) solutions of the
CL-Einstein-Dirac equation of type II. The essential idea of this
construction is the fact that, under conformal change of metrics,
the CL-Einstein-Dirac equation of type II behave in a relatively
stable way (more stable than the CL-Einstein-Dirac equation of
type I and both types of the classical Einstein-Dirac equation).
More precisely, we show in Section 4 that if $( Q^{n,r}, \eta )$
admits a non-trivial spinor field $\psi$, called a {\it reduced
weakly parallel spinor},  satisfying the differential equation in
Definition 4.3, then over the manifold $( Q^{n,r}, \overline{\eta}
= e^u \eta  )$ with conformally  changed metric $\overline{\eta} =
e^u \eta$ the pullback $\overline{\psi}$ of $\psi$ becomes a
CL-Einstein spinor of type II (see Theorem 4.2). Parallel spinors
[8] are trivial examples for reduced weakly parallel spinors. In
Section 5 we will provide examples for reduced weakly parallel
spinors that are not parallel spinors (see Theorem 5.2).

\bigskip   \noindent
\section{Coupling of the square of the Dirac operator to the Einstein equation}

\indent We first recall the process of obtaining the classical
Einstein-Dirac equation of type I in pseudo-Riemannian signature
[6, 7].  Applying the process to the behaviour of the square of
the Dirac operator under  change of metrics, we then derive the
classical Einstein-Dirac equation of type II.  \par Let $h$ be a
symmetric (0,2)-tensor field on $(Q^{n,r},  \eta)$, and let $H$ be
the (1,1)-tensor field induced by $h$ via $h (X,Y) = \eta ( X ,
H(Y) )$.
 Then  the tensor field  $\overline{\eta}$ defined by
\begin{equation}      \overline{\eta} (X , Y ) = \eta ( X ,  e^H(Y) ) = \eta( e^{\frac{H}{2}} (X) ,  e^{\frac{H}{2}} (Y) )    \end{equation}
is a pseudo-Riemannian metric of the same index $r$. Let $K  :=
e^{\frac{H}{2}} $ and  let  $ \Lambda$ be the (1,2)-tensor field
defined by
\begin{eqnarray*}
2 \, \eta ( \Lambda (X , Y ),  Z) &  =  &   \eta \left(  Z  , \, K
\{  ( \nabla^{\eta}_{K^{-1}(X)} K^{-1} ) (Y) \}  -
K \{  ( \nabla^{\eta}_{K^{-1}(Y)}  K^{-1} )(X) \}  \right)    \\
&     &     \\
&    & +  \eta \left(  Y  , \, K \{  ( \nabla^{\eta}_{K^{-1}(Z)}
K^{-1} )(X) \}  -
K \{  ( \nabla^{\eta}_{K^{-1}(X)}  K^{-1} )(Z) \}  \right)    \\
&     &     \\
&     & +  \eta \left(  X  , \, K \{  ( \nabla^{\eta}_{K^{-1}(Z)}
K^{-1} )(Y) \}  - K \{  ( \nabla^{\eta}_{K^{-1}(Y)}  K^{-1} )(Z)
\} \right).
\end{eqnarray*}
Then the Levi-Civita connections $\nabla^{\overline{\eta}}$ and
$\nabla^{\eta}$ are related by
\begin{equation}
\nabla^{\overline{\eta}}_{K^{-1}(X)} \left( K^{-1}(Y) \right) \ =
\ K^{-1} \left( \nabla^{\eta}_{K^{-1}(X)} Y \right) +  K^{-1}
\left\{ \Lambda (X, Y) \right\} .
\end{equation}
  Let $
{\widehat{K}} : \Sigma {(Q)}_{\overline{\eta}} \longrightarrow
\Sigma {(Q)}_{\eta} $ be a natural isomorphism  preserving the
inner product of spinors and the Clifford multiplication with
\begin{equation}
 \langle  \,  {\widehat{K}} ( \varphi)  , \, {\widehat{K}} ( \psi)   \,  \rangle_{\eta}   =   \langle  \varphi , \, \psi  \rangle_{\overline{\eta}} ,   \qquad
 ( K X ) \cdot ( {\widehat{K}} \psi )   =    {\widehat{K}} ( X \cdot \psi )      \end{equation}
for all $X \in \Gamma(T(Q))  , \  \varphi , \psi \in  \Gamma ( \Sigma {(Q)}_{\overline{\eta}} )$, where
 the dot "$\cdot$" in the latter relation indicates the Clifford multiplication with respect to $\eta$ and $\overline{\eta}$,
 respectively.
Let $( E_1 , \ldots , E_n )$ be a local $\eta$-orthonormal frame
field on $( Q^{ n , r} , \eta )$.    For shortness we introduce
the notation $\chi (i) : =  \eta (E_i, E_i)$ and $\chi ( i_1
\ldots i_s ) : = \chi(i_1) \, \chi(i_2) \cdots \chi (i_s)$ for $1
\leq s \leq n$. Then, because of (2.2), the spinor  derivatives
$\nabla^{\eta}, \, \nabla^{\overline{\eta}}$ are related by [4]
\begin{equation}
\Big\{  \widehat{K} \circ  \nabla^{\overline{\eta}}_{K^{-1}(E_j)}  \circ  \big( \widehat{K} \big)^{-1}   \Big\} (\psi)
  =    \nabla^{\eta}_{K^{-1}(E_j)} \psi  + \frac{1}{4} \sum_{k,l=1}^n \chi(k l) \Lambda_{jkl} E_k \cdot E_l \cdot \psi  ,
\end{equation}
where $\Lambda_{jkl} : = \eta ( \Lambda(E_j, E_k) , E_l )$, and
the Dirac operators $D_{\eta}, \, D_{\overline{\eta}}$  by
\begin{eqnarray}
&       &     \Big\{  \widehat{K} \circ D_{\overline{\eta}} \circ \big(\widehat{K} \big)^{-1}  \Big\} ( \psi )     \nonumber    \\
&       &      \nonumber      \\
&  =  &  \sum_{i=1}^n  \chi(i) E_i \cdot \nabla^{\eta}_{K^{-1}(E_i)} \psi + \frac{1}{4}
\sum_{j,k,l=1}^n  \chi(j k l) \Lambda_{jkl} E_j \cdot E_k \cdot E_l \cdot \psi         \nonumber      \\
&      &         \nonumber      \\
& = &  \sum_{i=1}^n  \chi(i) E_i \cdot \nabla^{\eta}_{K^{-1}(E_i)} \psi
- \frac{1}{2} \sum_{j,k=1}^n \chi(j k) \Lambda_{jjk} E_k \cdot \psi          \nonumber        \\
&     &         \nonumber        \\
&      &  + \frac{1}{2}  \sum_{j<k<l}^n  \chi(j k l) ( \Lambda_{jkl} + \Lambda_{klj} + \Lambda_{ljk} ) E_j \cdot E_k \cdot E_l \cdot \psi   .
\end{eqnarray}

\indent In order to  compute the infinitesimal variation of the
Dirac operator, we consider an one-parameter family of metrics of
index $r$,
\begin{equation}     \eta_t (X, Y)   : =  \eta ( X, e^{tH} (Y) )  = \eta ( e^{\frac{tH}{2}}(X) , e^{\frac{tH}{2}}(Y) ) ,    \qquad   \eta_o  := \eta,  \quad   t \in {\mathbb R}  ,     \end{equation}
which is generated by a symmetric (0,2)-tensor field $h$ on
$(Q^{n,r}, \eta)$. Let  $\Lambda_t$ be the (1,2)-tensor in (2.2)
determined by the pair
 $( \nabla^{\eta_t},   \nabla^{\eta})$ of the Levi-Civita connections  (with $K_t = e^{\frac{tH}{2}}$).  Let $\Omega_t$ be a 3-form generated by the tensor $\Lambda_t$ via
\begin{equation}     \Omega_t (X, Y, Z) = \eta ( \Lambda_t (X, Y), Z) + \eta ( \Lambda_t (Y, Z), X) + \eta ( \Lambda_t (Z, X), Y).     \end{equation}
Then direct computations show:

\begin{lem}
\begin{eqnarray*}
&     &  \frac{d}{dt} \,\bigg\vert_{t=0}  \Big\{ \Lambda_t (X, Y) - \Lambda_t (Y, X) \Big\}
 = - \frac{1}{2}  ( \nabla^{\eta}_X  H ) (Y) +  \frac{1}{2}  ( \nabla^{\eta}_Y  H ) (X) ,       \\
&     &     \\
&     &  \frac{d}{dt} \,\bigg\vert_{t=0}  \eta ( \Lambda_t (X, Y), Z )
=  \frac{1}{2}  \eta( ( \nabla^{\eta}_Y H ) (X), \, Z) -  \frac{1}{2}  \eta ( ( \nabla^{\eta}_Z H ) (X) , \, Y),     \\
&     &     \\
&     &   \frac{d}{dt} \,\bigg\vert_{t=0} \Omega_t (X, Y,  Z) = 0.
\end{eqnarray*}
\end{lem}

\noindent Applying  Lemma 2.1 to (2.5), we arrive at  the
variation formula of the Dirac operator:
\begin{eqnarray} &     &    \frac{d}{dt} \, \bigg \vert_{t=0} \,
\Big\{  \widehat{K_t} \circ D_{\eta_t} \circ \big(\widehat{K_t}
\big)^{-1}  \Big\} ( \psi )
\nonumber   \\
&     &      \nonumber      \\
& = &   -  \frac{1}{2} \sum_{j=1}^n \, \chi(j) \, h(E_j) \cdot \nabla^{\eta}_{E_j} \psi  -
 \frac{1}{4} {\rm div}_{\eta}(h) \cdot \psi +  \frac{1}{4} {\rm grad}_{\eta} ( {\rm Tr}_{\eta} (h) ) \cdot \psi  .
\end{eqnarray}

\bigskip  \noindent Recall [1] that for the standard complex product $\langle
\cdot , \cdot \rangle $ on the spinor bundle $\Sigma (Q)$, the
relation
\begin{equation}    \langle  X \cdot \varphi  ,  \psi    \rangle    +   (-1)^r    \langle   \varphi  ,  X \cdot \psi    \rangle     =    0       \end{equation}
holds for all vector fields $X$ and for all spinor fields $\varphi
, \psi$.  Taking the real part of (2.9) gives some simple but
crucial identities:
\begin{eqnarray}
(  ( \sqrt{-1})^r  X  \cdot \psi  ,  \psi  )   &   =   &    0   ,          \\
&        &      \nonumber      \\
( X \cdot \psi, Y \cdot \psi )   &  =    &    (-1)^r \eta(X, Y) (\psi, \psi) ,           \\
&       &       \nonumber       \\
( X \cdot Y \cdot \psi , \psi )   &   =   &   - \eta(X, Y ) (\psi,
\psi )  .
\end{eqnarray}

\noindent  Let {\rm Sym}(0,2) denote the space of all symmetric
(0,2)-tensor fields on $(Q^{n,r}, \eta)$, and let $(( \cdot ,
\cdot )) = (( \cdot , \cdot ))_{\eta}$  denote the naturally
induced metric on the space {\rm Sym}(0,2). Denote by
$\psi_{\eta_t} \, = \, \big( \widehat{K}_t \big)^{-1} (\psi)  \in
\Gamma ( \Sigma {(Q)}_{\eta_t} ) $ the pullback of $\psi =
\psi_{\eta} \, \in \Gamma ( \Sigma {(Q)}_{\eta} )$ via  natural
isomorphism $\widehat{K}_t$ (see (2.3)). Then (2.8) and (2.10)
together give the formula (1.4) for the first type energy-momentum
tensor $T_1$ :
\begin{equation}
 \frac{d}{dt} \,\bigg\vert_{t=0} \, \Big( \, {(\sqrt{-1})}^r  D_{\eta_t} \, \psi_{\eta_t} \, , \  \psi_{\eta_t} \, \Big)
\ = \ - \,  \frac{1}{4} \, (( \  T_1  \, , \ h \  ))   ,
\end{equation}
where
\begin{equation}
T_1  ( X , Y )   =   \Big( \, {(\sqrt{-1})}^r \{ X \cdot
\nabla_Y^{\eta} \psi + Y \cdot \nabla_X^{\eta} \psi \} \, , \ \psi
\, \Big).
\end{equation}

\noindent Moreover, using (2.8) and (2.9) and noting that
${(\sqrt{-1})}^r D_{\eta}$ is symmetric with respect to the
$L^2$-product, we can derive the formula (1.7) for the second type
energy-momentum tensor $T_2$.

\begin{lem} Let $U$ be an open subset of $Q^{n,r}$ with compact
closure, and let $h$ be a symmetric tensor field with support in
$U$. Then  for any spinor field $\psi$  on  $( Q^{n, r} , \eta )$,
we have
\[ \frac{d}{dt} \, \bigg \vert_{t=0} \,   \int_U  \Big( \,  (
D_{\eta_t}  \circ  D_{\eta_t} )  ( \psi_{\eta_t} ),  \,
\psi_{\eta_t} \Big) \, \mu_{\eta} =    -  \frac{1}{4} \, \int_U \,
(( \ T_2,   \  h  \  ))  \,  \mu_{\eta},
\]
where
\begin{eqnarray}
T_2  (X, Y)   &  =  &  \Big( X \cdot  \nabla^{\eta}_Y (D_{\eta}
\psi) +  Y \cdot  \nabla^{\eta}_X (D_{\eta}  \psi) ,   \    \psi
\Big)
\nonumber       \\
&         &       \nonumber       \\
&          &  + \,  (-1)^r   \Big( X \cdot  \nabla^{\eta}_Y \psi +
Y \cdot  \nabla^{\eta}_X  \psi ,   \   D_{\eta}  \psi  \Big) .
\end{eqnarray}
\end{lem}

\noindent {\bf Proof.}   Letting  $D = D_{\eta}$ and $\psi =
\psi_{\eta}$, we compute
\begin{eqnarray*}
&       &    \frac{d}{dt} \, \bigg \vert_{t=0} \,   \int_U  \Big( \,  ( D_{\eta_t}  \circ  D_{\eta_t} )  ( \psi_{\eta_t} ),  \,  \psi_{\eta_t} \Big)_{\eta_t} \,  \mu_{\eta}         \\
&       &        \\
&   =  &   \int_U  \,   \Big( \,    \frac{d}{dt}  \Big \vert_{t=0}
(\widehat{K}_t D_{\eta_t}) ( D \psi )_{\eta_t} ,   \   \psi  \,
\Big)    \mu_{\eta}
 +   \int_U   \,   \Big(    \,   D_{\eta}  \Big(   \frac{d}{dt}  \Big \vert_{t=0}  (\widehat{K}_t  D_{\eta_t})( \psi_{\eta_t})  \Big) ,   \    \psi  \, \Big)   \mu_{\eta}     \\
&       &        \\
&   =   &   \int_U    \Big(  -  \frac{1}{2} \sum_{j=1}^n \,
\chi(j) \, h(E_j) \cdot \nabla^{\eta}_{E_j} (D \psi)  -
 \frac{1}{4} {\rm div}_{\eta}(h) \cdot ( D \psi)  +  \frac{1}{4} {\rm grad}_{\eta} ( {\rm Tr}_{\eta} (h) ) \cdot (D \psi) ,   \
\psi   \Big)  \mu_{\eta}        \\
&        &        \\
&        &  +   \int_U    \Big(  ( \sqrt{-1} )^{3r}   \Big\{ -
\frac{1}{2} \sum_{j=1}^n \, \chi(j) \, h(E_j) \cdot
\nabla^{\eta}_{E_j}  \psi  -
 \frac{1}{4} {\rm div}_{\eta}(h) \cdot  \psi         \\
&         &         \\
&         &   \qquad   \quad   +  \frac{1}{4} {\rm grad}_{\eta} ( {\rm Tr}_{\eta} (h) ) \cdot  \psi    \Big\}  ,   \
( \sqrt{-1} )^r  D_{\eta}  \psi   \Big)   \mu_{\eta}       \\
&        &        \\
&   =    &   -   \frac{1}{2}   \,   \int_U     \Big(  \sum_{i=1}^n
\chi(i)  h(E_i)  \cdot   \nabla^{\eta}_{E_i}  (D  \psi) ,   \ \psi
 \Big)   \mu_{\eta}
   -   \frac{(-1)^r}{2}   \,   \int_U     \Big(  \sum_{i=1}^n   \chi(i)  h(E_i)   \cdot   \nabla^{\eta}_{E_i}   \psi ,   \   D  \psi
 \Big)   \mu_{\eta}     \\
&       &      \\
 &  =     &      -   \frac{1}{4}    \,    \int_U   \, (( \ T_2,   \  h  \  ))  \,  \mu_{\eta}  .
\end{eqnarray*}
\hfill{$\Box$}

\bigskip \noindent  We further need to recall the well-known
formulas for the variation of the volume form and the scalar
curvature, which one easily obtain from (2.6) and from the
pseudo-Riemannian version of the second formula in Proposition 2.2
of [7].

\begin{lem}  (see [3])  Let $U$ be an open subset of $Q^{n,r}$ with compact
closure, and let $h$ be a symmetric tensor field with support in
$U$.  Then  we have
\begin{eqnarray*}
 \frac{d}{dt} \, \bigg\vert_{t=0} \, \mu_{\eta_t} & = &  \frac{1}{2} \, (( \ \eta \, , \ h \ )) \, \mu_{\eta} ,   \\
&     &     \\
 \frac{d}{dt} \, \bigg\vert_{t=0} \,   \int_U  S_{\eta_t} \, \mu_{\eta}  &  = &  -
\int_U (( \, {\rm Ric}_{\eta}  , \ h \ )) \, \mu_{\eta}.
\end{eqnarray*}
\end{lem}

\noindent Making use of Lemma 2.2 and 2.3 and following the proof
of Theorem 2.1  of [6], we now establish the main result of this
section.

\begin{thm}       Let $Q^{n , r}$ be a pseudo-Riemannian spin manifold.
Fix the notation $P_{\eta}$ to mean either $P_{\eta} = ( \sqrt{-1}
)^r D_{\eta}$ or $P_{\eta} = D_{\eta}  \circ D_{\eta}$. Then,
  a pair $( \eta_o , \, \psi_o )$ is a critical point of the Lagrange functional
\[
  W ( \eta ,  \psi )
    =      \int_U   \Big\{ a  S_{\eta}   + b +  \epsilon    \nu     ( \psi_{\eta} ,  \psi_{\eta}  )_{\eta} -   \epsilon   (  P_{\eta} (\psi) , \, \psi )_{\eta}
\Big\}   \mu_{\eta}    ,    \qquad    a, b, \epsilon,  \nu     \in   {\mathbb R} ,     \quad   \epsilon  \neq 0 ,
\]
for all open subsets $U$ of $Q^{\, n,r}$ with compact closure  if
and only if $( \eta_o ,  \, \psi_o )$ is a solution of the
following system of differential equations:
\begin{equation}
  P_{\eta} ( \psi  )  =    \nu \,   \psi     \qquad
 \mbox{and}     \qquad
 a   \Big\{  {\rm Ric}_{\eta} -   \frac{1}{2}  \,   S_{\eta}  \eta  \Big\}  - \frac{b}{2}  \eta =  \frac{\epsilon}{4}  \,  T,
\end{equation}
where $T$ is a symmetric tensor field defined by (2.14) or by
(2.15) depending on a choice of $P_{\eta}$.
\end{thm}

\indent We close the section with generalizing Definition 2.1 and
3.1 of [6].

\bigskip \noindent {\bf Definition 2.1}  (i)  A non-trivial spinor field
$\psi$ on $(Q^{n,r} , \eta), \ n \geq 3 $, is called an {\it
Einstein spinor of type I}  for the  eigenvalue $( \sqrt{-1}
)^{3r}
\nu_1,  \, \nu_1  \in {\mathbb R},$   if it is a solution of the system (1.2)-(1.4).     \\
(ii)  A non-trivial spinor field $\psi$ on $(Q^{n,r} , \eta), \ n
\geq 3 $, is called an {\it Einstein spinor of type II}  for the
eigenvalue $\nu_2  \in  {\mathbb R}$   if it is a solution of the
system (1.6)-(1.7).

\bigskip \noindent
{\bf Definition 2.2}  Assume that  $a(n-2) S + bn \ (a, b  \in
{\mathbb R})$ does not vanish at any point of $( Q^{n,r}, \eta ),
\, n \geq 3$. A non-trivial spinor field $\psi$ on $( Q^{n,r},
\eta )$ is called a {\it weak Killing spinor}  (shortly,
WK-spinor) with WK-number  $( \sqrt{-1} )^{3r} \nu_1  \neq 0, \,
\nu_1  \in {\mathbb R},$ if $\psi$ is a solution of the
differential equation
\begin{equation}
\nabla_X  \psi    \  =  \  (\sqrt{-1})^{3r}  \beta(X)   \cdot
\psi   + n \, \alpha (X)  \psi   +  X  \cdot  \alpha  \cdot
\psi,
\end{equation}
where $\alpha$ is a  1-form and $\beta$ is a symmetric tensor
field defined by
\[
\alpha  =    \frac{a(n-2) \, dS}{2(n-1)  \{ a(n-2) S  + bn \} }
\]  and   \[ \beta     =      \frac{2 \, \nu_1}{ a(n-2) S + bn }
\Big\langle  a \Big\{ {\rm Ric}  -  \frac{1}{2}  S \eta \Big\} -
\frac{b}{2}  \eta   \Big\rangle ,
\]  respectively.

\bigskip  \noindent
{\bf Remark 2.1} As in the Riemannian case (see Theorem 3.1 of
[6]), any pseudo-Riemannian WK-spinor $\psi$  with positive length
$( \psi , \psi )
> 0$ (resp. negative length $( \psi ,  \psi ) < 0$) becomes an
Einstein spinor of type I: Since
\[   d \Big( \frac{(\psi, \psi )}{ a(n-2) S + bn }  \Big) = 0,     \]
it follows that
\[      \frac{(\psi, \psi )}{ a(n-2) S + bn }       \]
is constant on $Q^{n,r}$.
 One verifies easily that
 the equations (1.2)-(1.4) are indeed satisfied with
\[    \epsilon = - \, \frac{ a(n-2) S + bn }{\nu_1  \,  ( \psi ,  \psi )} .   \]

\bigskip   \noindent
{\bf   Remark 2.2}    Evidently, the solution space   of  the type
I classical Einstein-Dirac equation  is a subspace of that of the
type II classical Einstein-Dirac equation.  Hence it is of
interest to find such Einstein spinors of type II  that  are not
Einstein spinors of type I: Let $( Q^{n,r}, \eta)$ admit a spinor
field $\psi$ satisfying the differential equation [2]
\[     \nabla_X  \psi  = - (  \sqrt{-1} )^{3r+1} \,  \frac{\nu_1}{n} X  \cdot   \psi.     \]
Then the metric $\eta$ is necessarily  Einstein with scalar
curvature
\[    S =  (-1)^{r+1}  \,  \frac{4(n-1)  \nu_1^2 }{n} .      \]
If we choose the parameters $a$ and  $b$ so as to be related by
\[   b =  -  \frac{a(n-2)}{n} S   =  (-1)^r   \frac{4a(n-1)(n-2) \nu_1^2}{n^2},  \]
then  $\psi$ satisfies (1.6)-(1.7)  with
\[    \nu_2 =  (-1)^{r+1}  \,  \nu_1^2       \qquad       \mbox{and}      \qquad     a \Big\{ {\rm Ric}  -  \frac{S}{2}    \eta \Big\} -  \frac{b}{2}  \eta  =   \frac{\epsilon}{4}  T_2  =  0 .   \]
However, $\psi$ does not satisfy (1.2)-(1.4) in general.

\bigskip   \noindent
\section{Derivation of the  CL-Einstein-Dirac equations}

\indent Let  $\varphi= \varphi_{\eta}$ be a spinor field on
$(Q^{n,r}, \eta)$ such that  either $(\varphi, \varphi) > 0$ at
all points or $(\varphi, \varphi) < 0$ at all points.  We use the
simplifying notation
\[     \varphi^k  :=  ( \sigma  \varphi, \varphi )^k   \varphi,   \qquad   k \in {\mathbb R} ,    \]
where $\sigma = \sigma_{\varphi}  \in {\mathbb R}$ is a constant
defined by
\[    \sigma = 1 \   \mbox{if} \  (\varphi, \varphi) >0  \quad  \mbox{and}  \quad   \sigma = -1  \  \mbox{if} \    (\varphi, \varphi) < 0.     \]
Via direct computations, one verifies easily the following
variation formulas.

\noindent
\begin{lem} Let $U$ be an open subset of $( Q^{n, r} , \eta )$
with compact closure, and let $\varphi_c$ be a spinor field with
support in $U$. Then we have
\begin{eqnarray*}  &  (i)  &       \frac{d}{dt} \,\bigg\vert_{t=0}   ( \,  \sigma  \, ( \varphi + t \varphi_c )^k ,    \, ( \varphi + t \varphi_c )^k  \, )     =  2 (2k +1)   (\sigma  \varphi,  \varphi )^{2k} ( \sigma \varphi, \, \varphi_c ),      \\
&         &         \\
&  (ii)   &    \frac{d}{dt} \,\bigg\vert_{t=0}   \int_U    \Big(
\, \sigma  \, P_{\eta} ( \varphi + t \varphi_c )^k ,   \,    (
\varphi + t \varphi_c  )^k  \,
\Big)    \mu_{\eta}           \\
&          &        \\
&          &   =   4k   \int_U  \Big( \, \sigma \,   P_{\eta}   \{
( \sigma \varphi,  \varphi )^{k} \varphi   \},    \,
( \sigma \varphi, \varphi)^{k-1} \varphi  \Big)  ( \sigma \varphi,  \varphi_c )  \mu_{\eta}       \\
&          &         \\
&          &  \quad  +  2   \int_U   \Big( \, \sigma \,   P_{\eta}
\{ ( \sigma \varphi,  \varphi )^{k} \varphi   \},    \, ( \sigma
\varphi, \varphi)^{k}   \varphi_c     \Big)   \mu_{\eta},
\end{eqnarray*}
where $P_{\eta} = ( \sqrt{-1} )^r D_{\eta}$ or $P_{\eta} =
D_{\eta}  \circ  D_{\eta}$.
\end{lem}

\noindent
\begin{thm}
Let $Q^{n, r}$ be a pseudo-Riemannian spin manifold.  Consider the
Lagrange functional
\[
  W ( \eta ,  \varphi )
    =      \int_U   \Big\{ a  S_{\eta}   + b +  \epsilon    \nu     ( \sigma  \, \varphi^k ,  \varphi^k  )_{\eta} -   \epsilon   ( \sigma \, P_{\eta} (\varphi^k) , \, \varphi^k )_{\eta}
\Big\}   \mu_{\eta}     \] over open subsets $U$ of $Q^{n,r}$ with
compact closure,  where $ a, b, k , \epsilon,  \nu     \in
{\mathbb R},  \,  \epsilon \neq 0$, are real numbers.

\bigskip \noindent
(i)  In case of $2k +1 \not= 0$, a pair $(\eta^{\ast},
\varphi^{\ast})$ is a critical point of $W ( \eta, \varphi )$ for
all open subsets $U$ of $Q^{n,r}$ with compact closure if and only
if $(\eta^{\ast}, \varphi^{\ast} )$ is a solution of the following
system of differential equations:
\begin{equation}
 P_{\eta} ( \varphi^k  )  =   \nu \,   \varphi^k     \qquad
 \mbox{and}     \qquad
 a   \Big\{  {\rm Ric}_{\eta} -   \frac{1}{2}  \,   S_{\eta}  \eta  \Big\}  - \frac{b}{2}  \eta =  \frac{\epsilon}{4}  \,  T ,
\end{equation}
where $T$ is a symmetric tensor field  defined by
\begin{equation}
T  ( X , Y )  = T_1  ( X , Y )   =   \Big( \, \sigma \,
(\sqrt{-1})^r \{ X \cdot \nabla^{\eta}_Y  \varphi^k + Y \cdot
\nabla^{\eta}_X  \varphi^k \} \, , \  \varphi^k  \, \Big)
\end{equation}
 if $P_{\eta} = ( \sqrt{-1} )^r D_{\eta}$ and defined by
\begin{eqnarray}
T ( X , Y )  = T_2  (X, Y)   &  =  & \sigma  \Big(  X \cdot  \nabla^{\eta}_Y (D_{\eta} \varphi^k) +  Y \cdot  \nabla^{\eta}_X (D_{\eta}  \varphi^k) ,   \    \varphi^k  \Big)         \nonumber  \\
&         &      \nonumber     \\
&          &  + \,  \sigma  (-1)^r   \Big( X \cdot \nabla^{\eta}_Y
\varphi^k +  Y \cdot  \nabla^{\eta}_X  \varphi^k ,   \ D_{\eta}
\varphi^k \Big)
\end{eqnarray}
if $P_{\eta} = D_{\eta}  \circ D_{\eta}$, respectively.

\bigskip  \noindent
(ii) In case of $2k +1 = 0$,  a pair $(\eta^{\ast}, \varphi^{\ast}
)$ is a critical point of $W ( \eta, \varphi )$ for all open
subsets $U$ of $Q^{n,r}$ with compact closure if and only if
$(\eta^{\ast}, \varphi^{\ast})$ is a solution of the following
system of differential equations:
\begin{equation}
P_{\eta} ( \varphi^k  )  =   f \,   \varphi^k
\end{equation}
and
\begin{equation}
 a   \Big\{  {\rm Ric}_{\eta} -   \frac{1}{2}  \,   S_{\eta}  \eta  \Big\}  -  \frac{b + \epsilon \nu }{2}  \eta   =  \frac{\epsilon}{4}  \,  T -  \frac{\epsilon}{2} \, f \, \eta  ,
\end{equation}
where  $f : Q^{n,r}  \longrightarrow  {\mathbb R}$ is a
real-valued function and $T$  is a symmetric tensor field defined
by (3.2) or by (3.3) depending on a choice of $P_{\eta}$.
\end{thm}

\noindent {\bf Proof.}   Let $h$ be a symmetric tensor field with
support in $U$, and let $\varphi_c$ be a spinor field with support
in $U$. Let $\eta_t$ be an one-parameter family of metrics in
(2.6). Using Lemma 3.1, we compute at $t=0$ :
\begin{eqnarray*}
&    &   \frac{d}{dt}  W(\eta_t  ,\varphi + t  \varphi_c)  =  \frac{d}{dt}  W(\eta_t  , \varphi )  + \frac{d}{dt}  W(\eta  ,\varphi + t  \varphi_c )    \\
&     &        \\
&   = &  \frac{d}{dt}    \int_U   a  S_{\eta_t}  \mu_{\eta}   +       \frac{d}{dt}        \int_U   a  S_{\eta}  \mu_{\eta_t}
  +   \frac{d}{dt}   \int_U  b   \mu_{\eta_t}     +  \frac{d}{dt}  \int_U    \epsilon  \nu  ( \sigma \varphi^k , \varphi^k )  \mu_{\eta_t}       \\
&       &        \\
&       & - \frac{d}{dt} \int_U \epsilon ( \sigma P_{\eta} (
\varphi^k ), \, \varphi^k ) \mu_{\eta_t}
 - \frac{d}{dt}  \int_U  \epsilon ( \sigma P_{\eta_t} ( \varphi_{\eta_t}^k ),  \,  \varphi_{\eta_t}^k )    \mu_{\eta}      \\
&        &        \\
&        &     + \frac{d}{dt}  \int_U    \epsilon  \nu  ( \sigma  ( \varphi + t \varphi_c )^k ,  \, ( \varphi + t \varphi_c )^k )  \mu_{\eta}   -  \frac{d}{dt}  \int_U    \epsilon  ( \sigma  P_{\eta} ( \varphi + t \varphi_c )^k ,  \,  ( \varphi + t \varphi_c )^k )     \mu_{\eta}        \\
&        &        \\
&   =   &   \int_U    \Big( \Big(  - a {\rm Ric}_{\eta}  +
\frac{a}{2} S_{\eta}  \eta    + \frac{b}{2} \eta +
\frac{\epsilon}{4} T + \frac{\epsilon \nu}{2} ( \sigma \varphi^k ,
\varphi^k )  \eta
 -  \frac{\epsilon}{2} ( \sigma P_{\eta} ( \varphi^k ), \,  \varphi^k )   \eta,   \
 h  \,  \Big)   \Big)    \mu_{\eta}          \\
&       &         \\
&       &   +    \int_U    \Big(  2  \epsilon \nu (2k +1)  (  \sigma  \varphi ,  \varphi )^{2k}    \cdot  \sigma   \varphi  -   4 \epsilon k ( \sigma  \varphi ,  \varphi )^{-1}  ( \sigma \,  P_{\eta} (\varphi^k),  \, \varphi^k )  \cdot  \sigma   \varphi         \\
&       &       \\
 &      &     \qquad   \quad  - 2 \epsilon   (  \sigma  \varphi ,  \,   \varphi )^k   \cdot    \sigma   P_{\eta} ( \varphi^k ) ,    \  \varphi_c  \Big)
\mu_{\eta} .
 \end{eqnarray*}
It follows that a pair $(\eta^{\ast} ,  \varphi^{\ast} )$ is a
critical point of the
 functional $W(\eta ,  \varphi )$ for all open subsets $U$ of
$Q^{n,r}$ with compact closure  if and only if it is a solution of
the equations
\begin{equation}
 \frac{\epsilon}{4}  T =
  a {\rm Ric}_{\eta}   -   \frac{a}{2}  S_{\eta}  \eta
- \frac{b}{2}  \eta   - \frac{\epsilon \nu}{2} ( \sigma \varphi^k
, \varphi^k )  \eta +  \frac{\epsilon}{2} ( \sigma P_{\eta} (
\varphi^k ), \,  \varphi^k )   \eta
\end{equation}
and
\begin{equation}
P_{\eta} ( \varphi^k )    =   - 2k ( \sigma  \varphi ,   \varphi
)^{- 2k - 1} ( \sigma  P_{\eta} ( \varphi^k ),   \,  \varphi^k )
\, \varphi^k + \nu ( 2k +1 )  \,   \varphi^k .
\end{equation}
Inner product of (3.7) with $\sigma  \cdot \varphi^k$ gives
\begin{equation}
0 =   (2k +1)  \Big\{ ( \sigma   P_{\eta} ( \varphi^k ),   \,
\varphi^k ) -  \nu ( \sigma \varphi^k,  \varphi^k )  \Big\} ,
\end{equation}
and so, in case of $2k +1 \not= 0$, (3.6)-(3.8) imply  part (i) of
the theorem. Now we consider the other case  $2k +1 = 0$. In this
case, $( \sigma \varphi^k, \varphi^k ) = ( \sigma \varphi ,
\varphi )^{2k +1} = 1$  and hence (3.7) gives
\begin{equation}   P_{\eta} ( \varphi^k  )  =   f \,   \varphi^k   \end{equation}
with $f  := ( \sigma \, P_{\eta} (\varphi^k),  \, \varphi^k ) $.
Thus, (3.6) and (3.9) together prove  part (ii) of the theorem.
\hfill{$\Box$}

\bigskip  \noindent
We observe that the system (3.1)-(3.3) is  not new and  is in fact
equivalent to the classical system (2.16). We therefore focus our
attention on the system (3.4)-(3.5) which is a new Einstein-Dirac
system.

\bigskip   \noindent
{\bf  Definition 3.1}    A non-trivial spinor field $\psi$  on
$(Q^{n,r} , \eta), \ n \geq 3 $, is called a {\it CL-Einstein
spinor of type I } (resp. {\it type II}) {\it with characteristic
function} $f$ if it is of constant length $\vert \psi  \vert = \pm
1$ and satisfies the system (1.9) and (1.10) (resp.  (1.9) and
(1.11)).

\bigskip   \noindent
{\bf  Remark 3.1}  Let  $\varphi$ be a spinor field on $(Q^{n,r},
\eta)$ such that either $( \varphi, \varphi) > 0$ at all points or
$( \varphi, \varphi )  < 0$ at all points.  Let $T_1$ and $T_2$ be
symmetric tensor fields induced by $\varphi$ as in (1.10) and
(1.11), respectively. Then, via direct computations, one finds
that
\begin{eqnarray}
&       &  {\rm  div} (T_1)(X)    \ = \  \sigma   \sum_{i=1}^n  \chi(i)    ( \nabla_{E_i}  T_1 )( E_i, X)            \\
&       &      \nonumber    \\
&   =   &  \sigma  \Big( ( \sqrt{-1} )^r   \nabla_X (D \varphi ), \, \varphi  \Big)  -  \sigma  \Big(  \nabla_X \varphi,  \,  ( \sqrt{-1} )^r   D \varphi    \Big)
-  \sigma  \Big( ( \sqrt{-1} )^r  X  \cdot  D^2 \varphi , \, \varphi  \Big)       \nonumber
\end{eqnarray}
and
\begin{eqnarray}
  {\rm  div} (T_2)(X)    &  =   &
  \sigma   \Big(    \nabla_X (D^2 \varphi ), \, \varphi    \Big)  -  \sigma  \Big(  \nabla_X \varphi,  \,  D^2 \varphi      \Big)      \nonumber      \\
&       &        \nonumber      \\
&        &  -  \sigma   \Big(   X  \cdot  D^3 \varphi , \, \varphi
\Big)     -   (-1)^r  \sigma  \Big(   X  \cdot  D^2 \varphi , \, D
\varphi   \Big).
\end{eqnarray}

\noindent (i) If  $( \sqrt{-1} )^r   D \varphi  =  f_1  \varphi$
for some function $f_1 : Q^{n,r} \longrightarrow {\mathbb R}$ and
$\varphi$ is of constant length $\vert \varphi \vert =  \pm 1 $,
then
\[  {\rm div}(T_1)(X)
=    2  \,  d f_1(X) ( \sigma \varphi,  \varphi )   =  2  \, d
f_1(X) ,
\]
and so
\begin{equation}   {\rm div} \Big(  \frac{1}{4}  T_1 -  \frac{f_1}{2}  \eta \Big)  = 0 ,     \end{equation}
which is  required by the Einstein equation  in (1.9).     \\
\noindent   (ii) Similarly, if  $ D^2 \varphi  =  f_2  \varphi$
for some function $f_2 : Q^{n,r} \longrightarrow {\mathbb R}$ and
$\varphi$ is of constant length $\vert \varphi \vert =  \pm 1 $,
then
\[
  {\rm div}(T_2)(X)
=    2  \,  d f_2(X) ( \sigma \varphi,  \varphi )   =  2  \, d f_2
(X) ,
\]
and so
\begin{equation}   {\rm div} \Big(  \frac{1}{4}  T_2 -  \frac{f_2}{2}  \eta \Big)  = 0 .     \end{equation}
From  (3.12)-(3.13) we see that the Einstein equation  \[ \Big\{
{\rm Ric} - \frac{S}{2} \eta  \Big\}  -    \frac{c}{2} \eta =
\frac{\epsilon}{4}  T   - \frac{\epsilon}{2} \, f \, \eta  \] of
the CL-Einstein-Dirac equation (1.9) has a natural coupling
structure.  However, we should note  that neither (3.12) nor
(3.13) holds in general, unless $(\varphi, \varphi )$ is of
constant length.

\bigskip
We can rewrite the CL-Einstein-Dirac equation of type I
 \begin{eqnarray}          ( \sqrt{-1} )^r D  \psi  &  =  &  f_1
 \psi ,        \\
&       &     \nonumber   \\
     a  \Big\{ {\rm Ric}  -  \frac{S}{2}    \eta  \Big\}  -    \frac{c}{2} \eta  &
     =  &
 \frac{\epsilon}{4}  T_1   -  \frac{\epsilon}{2} \, f_1  \, \eta,
\end{eqnarray}
where
\begin{equation}
T_1 (X , Y) =  \Big( \sigma {(\sqrt{-1})}^r  \{ X \cdot \nabla_Y
\psi + Y \cdot \nabla_X \psi \}  , \, \psi  \Big) ,
\end{equation}
in an equivalent form:  Since contracting both sides of (3.15)
gives
\begin{equation}
\epsilon (n-1) f_1 = a (n-2) S + cn,
\end{equation}
one checks that the system (3.14)-(3.15) is actually equivalent to
the system
\begin{equation}
\epsilon  {(\sqrt{-1})}^r D \psi =    \Big\{   \frac{a(n-2)}{n-1} \, S  +  \frac{cn}{n-1}    \Big\}  \psi
\end{equation}
and
\begin{equation}
a  \Big\{ {\rm Ric}  -  \frac{S}{2(n-1)}    \eta  \Big\}   +
\frac{c}{2(n-1)} \eta =   \frac{\epsilon}{4}  T_1 .
\end{equation}

\noindent Since the system (3.18)-(3.19) is similar to the
classical Einstein-Dirac equation of type I, we are led to an
analogue of the WK-equation in Definition 2.2.

\bigskip   \noindent
{\bf  Definition 3.2}  A non-trivial spinor field $\psi$ on
$(Q^{n,r} , \eta), \ n \geq 3 $, is called a  {\it  WW-spinor}  if
$\psi$ satisfies the differential equation
\begin{equation}
\nabla_X  \psi  =   {(\sqrt{-1})}^{3r}  \Big(  - \frac{2a}{\epsilon}   \Big)  \Big\{   {\rm Ric}(X) -  \frac{S}{2(n-1)} X +
\frac{c}{2a(n-1)}  X   \Big\}  \cdot   \psi
\end{equation}
for some constants $\epsilon, a , c \in {\mathbb R} ,   \   \epsilon   \neq 0 ,    a  \neq  0 ,$ and for all vector fields $X$.

\bigskip  \noindent
Note that if the scalar  curvature $S$ of $(Q^{n,r} , \eta)$ is
constant, then the WW-equation (3.20) is  equivalent to the
WK-equation (2.17).  Because of (2.10), the length $\vert \psi
\vert$ of any WW-spinor $\psi$  is constant.  It follows that, by
rescaling the length $\vert \psi \vert$ if necessary, one may
assume without loss of generality that any WW-spinor $\psi$ is of
unit length $\vert \psi \vert = \pm 1$ or of zero length $\vert
\psi \vert = 0$.  As any WK-spinor of positive (resp. negative)
length is an Einstein spinor of type I, one then checks that any
WW-spinor $\psi$ of unit length is a CL-Einstein spinor of type I.

\bigskip   \noindent
\section{Constructing solutions of the CL-Einstein-Dirac equation of type II}

\indent     Let $\eta_1$ and $\eta_2$, $\eta_2 = e^u \eta_1$, be
conformally equivalent metrics on $Q^{n,r}$. By (2.3) there are
natural isomorphisms $j : T(Q) \longrightarrow T(Q)$ and $j :
\Sigma(Q)_{\eta_1} \longrightarrow \Sigma(Q)_{\eta_2}$ preserving
the inner products of vectors and spinors as well as the Clifford
multiplication:
\begin{eqnarray*}
&       &   \eta_2 (jX, jY) =  \eta_1 (X, Y),   \quad   \langle  j \varphi_1,  j \varphi_2  \rangle_{\eta_2} = \langle  \varphi_1,  \varphi_2  \rangle_{\eta_1} ,     \\
&       &     \\
&       &  (jX) \cdot (j \varphi) = j( X \cdot \varphi ),   \quad  X, Y  \in  \Gamma(T(Q)),   \quad    \varphi,  \varphi_1,  \varphi_2   \in  \Gamma(  \Sigma(Q)_{\eta_1} ).
\end{eqnarray*}
Denote by $\overline{X} := j(X)$ and $\overline{\varphi} := j (\varphi)$  the
corresponding vector fields and spinor fields on $(Q^{n,r}, \eta_2)$, respectively.  Then, for any spinor field $\psi$ on $(Q^{n,r}, \eta_1)$, we have
\begin{eqnarray}
    \nabla_{\overline{X}}^{\eta_2}  \,  \overline{\psi}
 &  =  &       e^{- \frac{u}{2}}  \,   \overline{ \nabla^{\eta_1}_X  \psi }  -  \frac{1}{4}
\eta_2 ( \overline{X},  \, {\rm grad}_{\eta_2} (u) ) \,  \overline{\psi}  -  \frac{1}{4}
\overline{X}  \cdot {\rm grad}_{\eta_2} (u)  \cdot  \overline{\psi} ,            \\
&      &     \nonumber       \\
D_{\eta_2}  \overline{\psi}   & =  &  e^{- \frac{u}{2}} \,  \overline{ D_{\eta_1} \psi }   +  \frac{n-1}{4} {\rm grad}_{\eta_2}(u)  \cdot
\overline{\psi} ,           \\
&      &      \nonumber        \\
( D_{\eta_2}   \circ  D_{\eta_2} )  \overline{\psi}  & =  &
e^{-u}  \,  \overline{ ( D_{\eta_1}  \circ D_{\eta_1} )  \psi }  -  \frac{1}{2} e^{- \frac{u}{2}}  {\rm grad}_{\eta_2}(u)  \cdot
\overline{ D_{\eta_1}  \psi}            \\
&        &      \nonumber         \\
&        &  -  \frac{n-1}{2}  e^{-u}  \,  \overline{  \nabla_{ {\rm grad}_{\eta_1}(u) } \psi }   +
\frac{ (n-1)^2 }{16}  \vert du  \vert_{\eta_2}^2  \,  \overline{\psi}  +  \frac{n-1}{4}  \triangle_{\eta_2} ( u)    \, \overline{\psi}.   \nonumber
\end{eqnarray}

\noindent Now consider a special class of spinors.

\bigskip   \noindent
{\bf Definition 4.1}    A non-trivial spinor field $\psi$ on
$(Q^{n,r} , \eta), \ n \geq 3 $, is called a {\it weakly
T-parallel spinor}  with confomal factor $u$  if it is of constant
length $\vert \psi  \vert = \pm 1$ and the equation
\begin{equation}      \nabla_X \psi \ = \  - \frac{1}{4}   du(X) \psi  - \frac{1}{4} \beta(X) \cdot  {\rm grad}(u)  \cdot  \psi
\end{equation}
holds for all vector fields $X$, for a symmetric (1,1)-tensor
field $\beta$ with
\[     {\rm Tr}( \beta ) = n,       \]
and for a real-valued function $u : Q^{n,r}  \longrightarrow
{\mathbb R}$ such that $\vert du \vert$ has no zeros on an open
dense subset of $Q^{n,r}$.

\bigskip  \noindent
Note that if $\psi$ is a parallel spinor on $(Q^{n,r}, \eta_1)$,
then the pullback $\overline{\psi}$ of $\psi$ is a weakly
T-parallel spinor on $(Q^{n,r}, \eta_2)$ with $\beta$ = the
identity map. In the following, we identify via the metric $\eta$
any exact 1-form "$du$" with the vector field "${\rm grad}(u)$"
and (1,1)-tensor field $\beta$ with the induced (0,2)-tensor field
$\beta(X, Y) = \eta(X, \beta(Y) )$.

\begin{pro}
Let $(Q^{n,r}, \eta)$ admit a weakly T-parallel spinor $\psi$
solving the equation (4.4).  Then we have

\bigskip  \noindent
(i) \ $ \beta(du) =   du$,

\bigskip  \noindent
(ii) \  $ \nabla_{du}   \psi  = 0 $,

\bigskip   \noindent
(iii)  \  $D \psi  =   \frac{n-1}{4}  \,  du  \cdot   \psi $,

\bigskip    \noindent
(iv)  \  $ D^2   \psi   =   \{  \frac{ (n-1)^2 }{16}  \vert du
\vert^2  +  \frac{n-1}{4}  \triangle u   \} \, \psi$,  \ where
$\triangle := - {\rm div} \circ {\rm grad}$,

\bigskip  \noindent
(v)  \  $S = \frac{1}{4}  \{  (n-1)^2 +1 -  \vert  \beta  \vert^2  \}  \vert  du  \vert^2  + (n-1) \triangle u $.
\end{pro}

\noindent
{\bf Proof.}  Since $( \sigma \psi,  \psi) =1$ is constant and $\beta$ is symmetric,
\[
0   =    \sigma ( \nabla_X  \psi,  \psi ) =   -  \frac{1}{4}  du(X)  +  \frac{1}{4}  \eta ( \beta(X),  {\rm grad}(u) )
  =    - \frac{1}{4} du(X)   +  \frac{1}{4}  \eta ( X,  \beta (du)
  ),
\]
which  proves part (i).   Using (ii)-(iii), we compute
\begin{eqnarray*}
D^2  \psi  & = &  \frac{n-1}{4}  D ( du  \cdot  \psi ) =   \frac{n-1}{4}  \triangle(u) \psi  - \frac{n-1}{2}  \nabla_{du}  \psi -  \frac{n-1}{4}  du   \cdot   D \psi    \\
&        &         \\
& =  &  \Big\{  \frac{n-1}{4}  \triangle u  +  \frac{ (n-1)^2
}{16}  \vert du  \vert^2    \Big\} \, \psi ,
\end{eqnarray*}
which proves part (iv). Substituting (iv) and (4.4)  into the
Schr\"{o}dinger-Lichnerowicz formula $D^2  \psi =  \triangle \psi
+  \frac{S}{4}  \psi$, one proves part (v).    \hfill{$\Box$}

\bigskip  \noindent
{\bf Remark 4.1} It is remarkable that when $ Q^{n,r}$ is a closed
manifold, the function $f_2 = \frac{ (n-1)^2 }{16}  \vert du
\vert^2  + \frac{n-1}{4} \triangle u$ in part (iv) of Proposition
4.1 cannot be constant: Suppose $f_2$ is a constant and hence an
eigenvalue of $D^2$. Then $f_2$ must be equal to a "positive"
constant $\lambda^2$ and
 for metric $\eta_1 := e^{-u} \eta$, we have $ \triangle_{\eta_1}
( u ) \, = \, \frac{n-3}{4} \vert du \vert_{\eta_1}^2 +
\frac{4}{n-1} \, \lambda^2  \, e^u . $ The last relation is
however a contradiction, since the left-hand side becomes zero
after integration.

\bigskip
Let $\psi$ be a  weakly  T-parallel spinor on $(Q^{n,r}, \eta)$
solving the equation (4.4). Then, a direct computation gives

\noindent
\begin{eqnarray*}
\frac{\epsilon}{4}  T_2 (X, Y)
&  =  &  \frac{\epsilon \sigma}{4}  \Big( X \cdot \nabla_Y (D \psi ) + Y \cdot  \nabla_X ( D \psi ),   \ \psi   \Big)       \\
&       &       \\
&       &  +   \frac{\epsilon \sigma}{4}  ( -1)^r  \Big( X \cdot  \nabla_Y  \psi +  Y \cdot  \nabla_X  \psi,    \   D \psi  \Big)     \\
&       &       \\
&   =   &  \frac{\epsilon \sigma (n-1)}{16}  \Big( X \cdot \nabla_Y ( du \cdot  \psi ) + Y \cdot  \nabla_X ( du  \cdot  \psi ),   \ \psi   \Big)       \\
&       &       \\
&       &  +   \frac{\epsilon \sigma (n-1)}{16}  ( -1)^r  \Big( X \cdot  \nabla_Y  \psi +  Y \cdot  \nabla_X  \psi,    \   du  \cdot  \psi  \Big)     \\
&       &      \\
&  =  &  \frac{\epsilon \sigma (n-1)}{16}  \Big(  X  \cdot  \nabla_Y du  \cdot \psi +  Y  \cdot  \nabla_X  du  \cdot \psi,    \   \psi  \Big)     \\
&      &      \\
&      & - \frac{\epsilon \sigma (n-1)}{64}  \Big(  X  \cdot du  \cdot  \Big\{  du(Y)  \psi + \beta(Y) \cdot  du \cdot \psi   \Big\}          \\
&      &       \\
&       &    \qquad \quad \qquad  + Y  \cdot du  \cdot  \Big\{  du(X)  \psi + \beta(X) \cdot  du \cdot \psi  \Big\},   \  \psi   \Big)        \\
&       &        \\
&       &    - \frac{\epsilon \sigma (n-1) }{64} (-1)^r   \Big(  X   \cdot  \Big\{  du(Y)  \psi + \beta(Y) \cdot  du \cdot \psi   \Big\}          \\
&      &       \\
&       &    \qquad  \quad \qquad  \qquad  + Y   \cdot  \Big\{  du(X)  \psi + \beta(X) \cdot  du \cdot \psi  \Big\},   \    du \cdot \psi   \Big)        \\
&      &        \\
&  =   &  -  \frac{\epsilon (n-1) }{8}   \eta ( X,  \nabla_Y du )   -  \frac{\epsilon (n-1)}{16}  du(X)  du(Y)   + \frac{\epsilon (n-1)}{16}  \vert  du  \vert^2
\beta (X, Y) .
\end{eqnarray*}

\bigskip \noindent
Guided by the last computation, one immediately proves:

\noindent
\begin{thm}
Let $\psi$ be a  weakly T-parallel spinor on $(Q^{n,r}, \eta)$
such that $\beta$ and $u$ are related to the Ricci tensor  and the
scalar curvature  of $( Q^{n,r}, \eta )$  by
\begin{eqnarray}
   \vert  du   \vert^2   \,   \beta (X, Y)
&  =  &   \frac{4}{n-2}   \Big\{  {\rm Ric}(X, Y) - \frac{1}{2}  S \, \eta(X, Y)   \Big\}  -  \frac{2c}{a(n-2)}  \eta(X, Y)    \nonumber     \\
&      &     \nonumber     \\
&      &  +   2 \, \eta (X, \, \nabla_Y (du)  )   +  du(X)  du(Y)   \nonumber     \\
&      &     \nonumber      \\
&       &  +    \Big\{   \frac{n-1}{2}  \vert du  \vert^2  + 2  \,
\triangle u   \Big\}  \,   \eta(X, Y) ,
\end{eqnarray}
where  \,  $ a, c  \in  {\mathbb R} ,   \,   a  \neq 0$, are real
numbers.  Then $\psi$ becomes a solution of  the CL-Einstein-Dirac
equation of type II (i.e., the system (1.9) and (1.11)), where the
characteristic function $f$ is given by
\[   f =  \frac{ (n-1)^2 }{16}  \vert du  \vert^2  +  \frac{n-1}{4}  \triangle u    \]
and  the parameter $\epsilon$   should be chosen to satisfy
\[    \epsilon =   \frac{4 a(n-2)}{n-1} .    \]
\end{thm}

\bigskip   \noindent
{\bf  Definition 4.2}    A non-trivial spinor field $\psi$ on
$(Q^{n,r} , \eta), \ n \geq 3 $, is called a {\it  weakly parallel
spinor} (shortly, WP-spinor) with conformal factor $u$  if it is a
weakly  T-parallel spinor  with conformal factor $u$  and
satisfies (4.5) for some constants $ a , c  \in  {\mathbb R}, \, a
\neq 0$.

\bigskip   \noindent
{\bf  Definition 4.3}    A non-trivial spinor field $\psi$ on
$(Q^{n,r} , \eta), \ n \geq 3 $, is called a {\it reduced weakly
parallel spinor} (shortly, reduced   WP-spinor)  with  conformal
factor $u$  if it is of constant length $\vert  \psi  \vert  = \pm
1$ and the differential equation
\begin{equation}
\vert   du  \vert^2   \,  \nabla_X  \psi   \ = \ - \frac{1}{n-2}
\Big\{    {\rm  Ric} (X) -   \frac{S}{n}    X \Big\}   \cdot  du
\cdot   \psi
\end{equation}
holds for all vector fields $X$ and  for a real-valued  function
$u : Q^{n,r} \longrightarrow  {\mathbb R}$ with such properties
that $\vert du \vert$ has no zeros on an open dense subset of
$Q^{n,r}$ and $e^u$ is proportional  to the scalar curvature $S$,
i.e.,
\begin{equation}  S = c^{\ast} e^u,   \qquad    c^{\ast}   \in  {\mathbb R}.
\end{equation}

\bigskip  \noindent
Note that (4.6) generalizes the equation $\nabla_X \psi = 0$ for
parallel spinors and that any reduced WP-spinor $\psi$ is a
harmonic spinor $D \psi = 0$.  Applying (4.6) to $0 = \sigma \cdot
\vert du \vert^2 (  \nabla_X  \psi ,  \psi )$, one shows:

\noindent
\begin{pro}
Let $( Q^{n,r},  \eta )$ admit a  reduced  WP-spinor  $\psi$  with
conformal factor $u$.  Then
\[     \nabla_{ du }  \psi  = 0    \qquad      \mbox{and}      \qquad      {\rm Ric} (du) =  \frac{S}{n}  du .      \]
\end{pro}

\indent We are going to prove that the equation (4.5) for
WP-spinors is conformally equivalent to the equation (4.6) for
reduced WP-spinors.  Consider  conformally equivalent metrics
$\eta_2  = e^u \eta_1$ on $Q^{n,r}$. Let $(F_1, \ldots, F_n )$ be
a local $\eta_1$-orthonormal frame field on $Q^{n,r}$. Then $(
\overline{F}_1 := e^{- \frac{u}{2}} F_1,  \ldots, \overline{F}_n
:= e^{- \frac{u}{2}} F_n )$ is $\eta_2$-orthonormal. Since the
Ricci tensors ${\rm Ric}_{\eta_2}$ and ${\rm Ric}_{\eta_1}$ are
related by
\begin{eqnarray*}
&      & {\rm Ric}_{\eta_2} ( \overline{F}_i, \overline{F}_j ) -  e^{-u}  {\rm Ric}_{\eta_1}  (F_i, F_j )   \\
&      &       \\
&  =  & - \frac{n-2}{2}  \,  \eta_2 ( \overline{F}_i, \,  \nabla^{\eta_2}_{\overline{F}_j} ( {\rm grad}_{\eta_2}  u ) )   -   \frac{n-2}{4}  \,  du ( \overline{F}_i) du ( \overline{F}_j)      \\
&      &        \\
&      &   + \frac{1}{2}  \,  \triangle_{\eta_2}(u) \,  \eta_2 ( \overline{F}_i,  \overline{F}_j)   +  \frac{n-2}{4}   \,  \vert du  \vert_{\eta_2}^2  \,   \eta_2 ( \overline{F}_i,   \overline{F}_j)
\end{eqnarray*}
and the scalar curvatures $S_{\eta_2}$ and $S_{\eta_1}$ by
\[   S_{\eta_2} -  e^{-u} S_{\eta_1}          =   (n-1) \,
\triangle_{\eta_2} (u)  + \frac{(n-1)(n-2)}{4}  \vert du
\vert_{\eta_2}^2  ,
\]
we have in particular the following formula.

\noindent
\begin{lem}
\begin{eqnarray*}
&        &  {\rm Ric}_{\eta_2} ( \overline{F}_i, \overline{F}_j ) -  \frac{1}{2}  S_{\eta_2}  \,  \eta_2  ( \overline{F}_i,   \overline{F}_j)     \\
&         &        \\
&   =   &  e^{-u}   \Big\{   {\rm Ric}_{\eta_1}  (F_i, F_j )   -  \frac{1}{2}   S_{\eta_1}  \eta_1 (F_i, F_j)     \Big\}        \\
&         &        \\
 &         & - \frac{n-2}{2}  \,  \eta_2 ( \overline{F}_i, \,  \nabla^{\eta_2}_{\overline{F}_j} ( {\rm grad}_{\eta_2} u ) )   -   \frac{n-2}{4}  \,  du( \overline{F}_i) du ( \overline{F}_j)      \\
&      &        \\
&      &   -  \frac{n-2}{2}  \,  \triangle_{\eta_2}(u) \,  \eta_2 ( \overline{F}_i,  \overline{F}_j)   -  \frac{(n-2)(n-3)}{8}   \,  \vert du  \vert_{\eta_2}^2  \,   \eta_2  ( \overline{F}_i,   \overline{F}_j)  .
\end{eqnarray*}
\end{lem}

\noindent
\begin{thm}
   A non-trivial spinor field $\psi$ on $(Q^{n,r} , \eta_1)$ is a reduced WP-spinor with conformal factor $u$ if and only if
the pullback $\overline{\psi}$ of $\psi$  is a   WP-spinor on $(
Q^{n,r}, \eta_2 = e^u \eta_1 )$ with conformal  factor $u$.
\end{thm}

\noindent {\bf Proof.}   We first prove the necessity. Let  $\psi$
be a reduced  WP-spinor  on $(Q^{n,r}, \eta_1)$ with conformal
factor $u$. In the notation of  (4.1), we have
\begin{eqnarray*}
&       &   \vert du  \vert_{\eta_2}^2   \,   \nabla_{\overline{X}}^{\eta_2}  \,  \overline{\psi}        \\
&       &      \\
&   =   &  -  \frac{1}{n-2} \, e^{-u}  \Big\{  \overline{ {\rm Ric}_{\eta_1}(X) } -  \frac{1}{n} S_{\eta_1}  \, \overline{X}  \Big\}
\cdot  {\rm grad}_{\eta_2}(u)  \cdot   \overline{\psi}      \\
&       &      \\
&       &  -  \frac{1}{4}  \vert du  \vert_{\eta_2}^2
\eta_2 ( \overline{X},  \, {\rm grad}_{\eta_2} (u) ) \,  \overline{\psi}  -  \frac{1}{4}  \vert du  \vert_{\eta_2}^2  \,
\overline{X}  \cdot {\rm grad}_{\eta_2} (u)  \cdot  \overline{\psi}
\end{eqnarray*}
and hence
\begin{equation}
\nabla_{\overline{X}}^{\eta_2}  \,  \overline{\psi}    =  -  \frac{1}{4}  \eta_2 ( \overline{X},  \, {\rm grad}_{\eta_2} (u) ) \,  \overline{\psi}
-  \frac{1}{4}  \gamma( \overline{X} )   \cdot {\rm grad}_{\eta_2} (u)  \cdot  \overline{\psi} ,
\end{equation}
where $\gamma$ is a symmetric tensor field defined by
\begin{eqnarray}
&        &    \vert  d u \vert_{\eta_2}^2  \,  \gamma ( \overline{X}, \, \overline{Y} )    \nonumber     \\
&       &     \nonumber     \\
&   =   &   \frac{4}{n-2}  \, e^{-u}   \Big\{  {\rm Ric}_{\eta_1} (X, Y)  -  \frac{1}{n} S_{\eta_1}  \eta_1 (X, Y)  \Big\}
+   \vert  d u \vert_{\eta_2}^2  \,  \eta_1 (X, Y) .
\end{eqnarray}
On the other hand, using Lemma 4.1, we compute
\begin{eqnarray*}
\Phi ( \overline{X}, \, \overline{Y} )  &  :=  &
 \frac{4}{n-2}   \Big\{  {\rm Ric}_{\eta_2} ( \overline{X}, \, \overline{Y} ) - \frac{1}{2}  S_{\eta_2}  \, \eta_2 ( \overline{X}, \,  \overline{Y} )   \Big\}  -  \frac{2c}{a(n-2)}  \eta_2 ( \overline{X},  \,  \overline{Y} )         \\
&      &          \\
&      &  +   2 \, \eta_2 ( \overline{X}, \, \nabla^{\eta_2}_{\overline{Y}} ( {\rm grad}_{\eta_2} (u) )  )   +  du( \overline{X})  du( \overline{Y})       \\
&      &           \\
&       &  +    \Big\{   \frac{n-1}{2}  \vert du  \vert_{\eta_2}^2  + 2  \,  \triangle_{\eta_2} ( u)   \Big\}  \,   \eta_2 ( \overline{X}, \,  \overline{Y})     \\
&  =   &    \frac{4  \,  e^{-u} }{n-2}   \Big\{  {\rm
Ric}_{\eta_1} ( X, Y) - \frac{1}{2}  S_{\eta_1}   \eta_1 ( X, Y )
-  \frac{c \, e^u}{2a}  \eta_1(X, Y) \Big\}   + \vert du
\vert_{\eta_2}^2   \eta_1 (X, Y) .
\end{eqnarray*}
Choose the parameters $a, c \in {\mathbb R}$ such that the
constant $c^{\ast}$ in (4.7) satisfies
\[    c^{\ast} =  - \frac{cn}{a(n-2)} .     \]
Then $S_{\eta_1} = -  \frac{cn}{a(n-2)} \, e^u$ and
\begin{equation}
\Phi ( \overline{X}, \, \overline{Y} )    =    \vert  du
\vert_{\eta_2}^2  \,  \gamma ( \overline{X}, \, \overline{Y} ) .
 \end{equation}
From (4.8)-(4.10), we conclude that  $\overline{\psi}$ is a weakly
T-parallel spinor on   $( Q^{n,r},  \eta_2 = e^u \eta_1 )$
satisfying (4.5), i.e., $\overline{\psi}$ is a   WP-spinor. In
order to prove the sufficiency, we reverse the process of the
proof  for the necessity: Let $\overline{\psi}$  be a WP-spinor on
$( Q^{n,r},  \eta_2 = e^u \eta_1 )$.     Then we have
\begin{eqnarray*}
&        &  \vert  du  \vert_{\eta_2}^2  \,  \beta ( \overline{X}, \, \overline{Y} )          \\
&        &         \\
&   =    &   \frac{4  \,  e^{-u} }{n-2}   \Big\{  {\rm Ric}_{\eta_1} ( X, Y) - \frac{1}{2}  S_{\eta_1}   \eta_1 ( X, Y )  -  \frac{c \, e^u}{2a}  \eta_1(X, Y) \Big\}   + \vert du  \vert_{\eta_2}^2   \eta_1 (X, Y) .
\end{eqnarray*}
Contracting  both sides of this equation gives
\[     S_{\eta_1}  =  -  \frac{cn}{a(n-2)}  \, e^u .     \]
Using (4.1), one verifies that $\psi$  satisfies the equation
(4.6) indeed.
 \hfill{$\Box$}

\bigskip   \noindent
\section{An existence theorem  for  WK-spinors and that for  reduced WP-spinors}

\indent We show that every parallel spinor may evolve to a
WK-spinor (resp. a reduced WP-spinor).  We give a description for
the evolution in a more general way than that given in Section 5
of [7].  \par  Let $(M^n, g_M )$ be a Riemannian manifold, and let
$( {\mathbb R}, g_{\mathbb R} )$ be the real line with the
standard metric. Let $( Q^{n+1} = M^n \times {\mathbb R}, \,
\eta_1 = g_M + \chi (n+1)  g_{\mathbb R}  ), \, \chi(n+1) =  \pm
1,$ be the pseudo-Riemannian product manifold.  We will write
$g_{\mathbb R} = dt \otimes dt$ using the standard coordinate $t
\in {\mathbb R}$ and regard $\eta_1$ as a reference metric on
$Q^{n+1}$. Let $(F_1, \ldots, F_n)$  denote a local
$\eta_1$-orthonormal frame field on $(M^n, g_M)$ as well as its
lift to $(Q^{n+1}, \eta_1)$. Let $F_{n+1} = \frac{d}{dt}$ denote
the unit vector field on $({\mathbb R}, g_{\mathbb R})$ as well as
the lift to $(Q^{n+1}, \eta_1)$.  We consider a  doubly warped
product of $g_M$ and $g_{\mathbb R}$:
\begin{equation}    \eta_2 =  A^2 \Big( \sum_{i=1}^n F^i \otimes
F^i \Big) +  \chi(n+1)  B^2 dt \otimes dt ,  \end{equation} where
 $A = A(t), \, B = B(t) : {\mathbb R} \longrightarrow {\mathbb
R}$ are positive  functions on ${\mathbb R}$ and $\{ F^i =  \eta_1
( F_i,  \, \cdot ) \}$ is the dual frame field of $\{ F_i  \}$.
Let $g_{M_t}$ be the metric on slice $M_t := M^n \times \{ t \},
\, t  \in {\mathbb R},$ of the foliation $( Q^{n+1} = M^n \times
{\mathbb R}, \, \eta_1)$ induced by the reference metric $\eta_1$,
and let
 $\nabla^{g_{M_t}}$  be the Levi-Civita connection.  Then the Levi-Civita connection $\nabla^{\eta_2}$ of $(Q^{n+1}, \eta_2)$ is related to
$\nabla^{g_{M_t}}$ by
\begin{eqnarray}
\nabla^{\eta_2}_{\overline{F}_i}  \overline{F}_j     &  =  &   A^{-2}  \nabla^{g_{M_t}}_{F_i} F_j   -  \chi(n+1)  \,   \delta_{ij} \,  B^{-2}  A^{-1} A_t \,  F_{n+1},      \\
&         &     \nonumber      \\
\nabla^{\eta_2}_{\overline{F}_{n+1}}  \overline{F}_j    & = &
\nabla^{\eta_2}_{\overline{F}_{n+1}}  \overline{F}_{n+1}  \ = \ 0,
\qquad  \qquad   1 \leq i,j \leq n,
\end{eqnarray}
where  $( \overline{F}_1:=  A^{-1} F_1, \ldots, \overline{F}_n :=
A^{-1} F_n, \overline{F}_{n+1} :=  B^{-1} F_{n+1} )$ is a
$\eta_2$-orthonormal frame field and $A_t$ indicates the
derivative $A_t = dA (F_{n+1})$. The second fundamental form
$\Theta_{\eta_2} = - \nabla^{\eta_2} \overline{F}_{n+1}$ of slice
 $M_t$ is expressed as
\begin{equation}
 \Theta_{\eta_2} ( \overline{F}_j )  = -  B^{-1} A^{-1} A_t  \,  \overline{F}_j ,    \qquad    1 \leq j  \leq  n.
\end{equation}
Furthermore, the Ricci tensor ${\rm Ric}_{\eta_2}$ and the scalar
curvature $S_{\eta_2}$ of $(Q^{n+1}, \eta_2)$ are related to the
Ricci tensor ${\rm Ric}_{g_{M_t}}$ and the scalar curvature
$S_{M_t}$ of  slice $( M_t,  g_{M_t})$ by
\begin{eqnarray}
{\rm Ric}_{\eta_2} ( \overline{F}_i, \overline{F}_j )  &  =  &
A^{-2}  \,  {\rm Ric}_{g_{M_t}} ( F_i, F_j )  -  \chi(n+1) \,
(n-1) B^{-2} A^{-2}
A_t A_t  \delta_{ij}        \nonumber   \\
&        &     \nonumber     \\
&        &  +   \chi(n+1)  \{   B^{-3} A^{-1} B_t A_t -   B^{-2} A^{-1}  A_{tt}    \}  \delta_{ij},     \\
&        &      \nonumber    \\
{\rm Ric}_{\eta_2} ( \overline{F}_{n+1}, \overline{F}_{n+1} )   &  =  & n  B^{-2} A^{-1} (  B^{-1} B_t A_t - A_{tt} ),       \\
&        &       \nonumber     \\
{\rm Ric}_{\eta_2} ( \overline{F}_i, \overline{F}_{n+1} )   &  =  &  0,       \\
&       &      \nonumber      \\
S_{\eta_2}   &  =  &  A^{-2}  \, S_{g_{M_t}}  -   \chi(n+1)  \,  n(n-1)   B^{-2} A^{-2} A_t A_t      \nonumber      \\
&        &      \nonumber      \\
&        &   +  \chi(n+1)  \,  2n   \{  B^{-3} A^{-1} B_t A_t -
 B^{-2} A^{-1} A_{tt}  \},
\end{eqnarray}
where $A_{tt} = (A_t)_t$ indicates the second derivative.  From
now on,  we are  interested in a special case that the warping
functions $A$ and $B$ are related by
\begin{equation}
B = ( A^p )_t = p A^{p-1} A_t,  \qquad   p \neq 0  \in {\mathbb
R}.
\end{equation}

\bigskip  \noindent
{\bf Definition 5.1}  A doubly warped product  (5.1) is called a
{\it (Y)-warped product } of $(M^n, g_M)$ and $( {\mathbb R},
g_{\mathbb R})$ with warping function $A$ and (Y)-factor $p$ if
the relation (5.9) is satisfied for some constant $p \neq 0 \in
{\mathbb R}$.

\noindent
\begin{pro}
Let $(Q^{n+1}= M^n  \times  {\mathbb R}, \eta_2)$ be a (Y)-warped
product  of $(M^n, g_M)$ and $( {\mathbb R}, g_{\mathbb R} )$ with
warping function $A$ and (Y)-factor $p$. Then the formulas
(5.4)-(5.8) simplify to

\bigskip  \noindent
(i) \ $\Theta_{\eta_2} ( \overline{F}_i,  \overline{F}_j )  = -
p^{-1} A^{-p}  \delta_{ij},   \qquad  1 \leq i,j \leq n,$

\bigskip  \noindent
(ii) \  ${\rm Ric}_{\eta_2} ( \overline{F}_i, \overline{F}_j )   =
A^{-2}  \,  {\rm Ric}_{g_{M_t}} ( F_i, F_j ) +  \chi(n+1) \, (p-n)
p^{-2}  A^{-2p}   \delta_{ij}$,

\bigskip  \noindent
(iii)  \ ${\rm Ric}_{\eta_2} ( \overline{F}_{n+1},
\overline{F}_{n+1} )  = n (p-1)  p^{-2} A^{-2p}$,

\bigskip   \noindent
(iv)  \ ${\rm Ric}_{\eta_2} ( \overline{F}_i, \overline{F}_{n+1} )   = 0$,

\bigskip  \noindent
(v)  \  $S_{\eta_2}     =    A^{-2}  \, S_{g_{M_t}}   + \chi(n+1)
\,  n(2p-n-1)  p^{-2}  A^{-2p}$.
\end{pro}

\bigskip  \noindent
An argument similar to that of  Proposition 5.1 of [7] shows:

\noindent
\begin{pro}
Let $(Q^{n+1}= M^n  \times  {\mathbb R}, \eta_2)$ be a (Y)-warped
product  of $(M^n, g_M)$ and $( {\mathbb R}, g_{\mathbb R} )$ with
warping function $A$ and  (Y)-factor $ \frac{n}{2} $. Assume that
$(M^n,  g_M)$  is Ricci-flat. Then the weak Killing equation
(2.17), in case of $b=0$, is equivalent to the system of
differential equations
\[   \nabla^{g_{M_t}}_V  \psi  = 0  \qquad   \mbox{and}   \qquad
\nabla^{\eta_2}_{\overline{F}_{n+1}}  \psi  =  - ( \sqrt{-1} )^{3r}  \nu_1 \overline{F}_{n+1}  \cdot   \psi  +  \frac{1}{2}
{\rm Tr}_{g_{M_t}} ( \Theta_{\eta_2} )  \psi,      \]
where $V$ is an arbitrary vector field on $Q^{n+1}$ with $\eta_2 (V, \overline{F}_{n+1} ) = 0$.
\end{pro}

\noindent
\begin{pro}
Let $(Q^{n+1}= M^n  \times  {\mathbb R}, \eta_2)$ be a  (Y)-warped
product  of $(M^n, g_M)$ and $( {\mathbb R}, g_{\mathbb R} )$ with
warping function $A$ and (Y)-factor $ \frac{n+1}{2} $. Assume that
$(M^n,  g_M)$  is Ricci-flat. Then the reduced WP-equation in
Definition 4.3 (in case that we set $u= - \log A$) is equivalent
to the system of differential equations
\[   \nabla^{g_{M_t}}_V  \psi  = 0  \qquad   \mbox{and}   \qquad
\nabla^{\eta_2}_{\overline{F}_{n+1}}  \psi  =    \frac{1}{2}
{\rm Tr}_{g_{M_t}} ( \Theta_{\eta_2} )  \psi,      \]
where $V$ is an arbitrary vector field on $Q^{n+1}$ with $\eta_2 (V, \overline{F}_{n+1} ) = 0$.
\end{pro}

\noindent {\bf Proof.}  Since $u= - \log A$, we have
\begin{eqnarray*}    \vert  du   \vert_{\eta_2}^2    &  =   &  \chi(n+1)  \,  p^{-2}   A^{-2p}  ,   \\
&        &       \\
{\rm grad}_{\eta_2} (u)   &  =   &   -  \chi(n+1)   \,  p^{-1}
A^{-p}  \,  \overline{F}_{n+1} .
\end{eqnarray*}
Moreover, by part (v) of Proposition 5.1, the scalar curvature
$S_{\eta_2} =0$ vanishes. Thus the reduced WP-equation becomes
\begin{eqnarray}
\nabla^{\eta_2}_V  \psi     &  =     &  - \frac{1}{n-1}  {\rm
Ric}_{\eta_2} (V)
\cdot    \frac{ {\rm grad}_{\eta_2}(u) }{    \vert  du   \vert_{\eta_2}^2 }   \cdot   \psi     \nonumber    \\
&       &   \nonumber     \\
&  =   &   \frac{p}{n-1}   A^p  \,  {\rm Ric}_{\eta_2} (V)   \cdot   \overline{F}_{n+1}   \cdot    \psi     \nonumber    \\
&      &      \nonumber      \\
&  =   & -  \chi(n+1)  \,    \frac{1}{n+1}   A^{- \frac{n+1}{2}}
\,  V \cdot  \overline{F}_{n+1}  \cdot   \psi
\end{eqnarray}
and
\begin{eqnarray}
\nabla^{\eta_2}_{\overline{F}_{n+1}}  \psi     &   =    &   -
\frac{1}{n-1}  {\rm Ric}_{\eta_2} ( \overline{F}_{n+1})
\cdot    \frac{ {\rm grad}_{\eta_2}(u) }{    \vert  du   \vert_{\eta_2}^2 }   \cdot   \psi     \nonumber    \\
&       &   \nonumber     \\
&  =   &  -  \frac{n}{n+1}  A^{- \frac{n+1}{2}}   \psi   \ =  \
\frac{1}{2} {\rm Tr}_{g_{M_t}} ( \Theta_{\eta_2} )  \psi.
\end{eqnarray}
On the other hand,
\begin{eqnarray}
\nabla^{\eta_2}_V  \psi     &   =    &  \nabla^{g_{M_t}}_V  \psi +
\chi(n+1)  \,  \frac{1}{2}  \Theta_{\eta_2} (V)
\cdot  \overline{F}_{n+1} \cdot   \psi       \nonumber    \\
&        &       \nonumber      \\
&  =  &   \nabla^{g_{M_t}}_V  \psi   -  \chi(n+1)  \,
\frac{1}{n+1} A^{- \frac{n+1}{2}}  V  \cdot \overline{F}_{n+1}
\cdot   \psi .
\end{eqnarray}
From (5.10)-(5.12) we conclude the proof. \hfill{$\Box$}

\bigskip \noindent Following a standard argument in the proof of Proposition 5.2 and
Theorem 5.1 of [7] in pseudo-Riemannian signature, we now
establish the following existence theorems.

\noindent
\begin{thm}
Let $(Q^{n+1}= M^n  \times  {\mathbb R}, \eta_2)$ be a (Y)-warped
product  of $(M^n, g_M)$ and $( {\mathbb R}, g_{\mathbb R} )$ with
(Y)-factor $ \frac{n}{2}$.  If $(M^n, g_M)$ admits a parallel
spinor, then for any real number $\lambda_Q \in {\mathbb R} \neq 0
$, $(Q^{n+1}, \eta_2)$ admits a WK-spinor to WK-number $ (
\sqrt{-1} )^{3r} \lambda_Q$, where $r=0$ if $\chi(n+1) = 1$ and
$r=1$ if $\chi(n+1) = -1$, respectively.
\end{thm}

\noindent
\begin{thm}
Let $(Q^{n+1}= M^n  \times  {\mathbb R}, \eta_2)$ be a (Y)-warped
product  of $(M^n, g_M)$ and $( {\mathbb R}, g_{\mathbb R} )$ with
(Y)-factor $ \frac{n+1}{2}$. If $(M^n, g_M)$ admits a parallel
spinor, then  $(Q^{n+1}, \eta_2)$ admits a reduced WP-spinor that
is not a parallel spinor.
\end{thm}

\bigskip  \noindent
Theorem 5.1 above improves  Theorem 5.1 of [7], since  (Y)-warped
products  of $(M^n, g_M)$ and $( {\mathbb R}, g_{\mathbb R} )$
with (Y)-factor $\frac{n}{2}$ essentially generalize the metrics
in Lemma 5.3 of [7].

\bigskip

\begin{Literature}{xx}
\bibitem{1}
H. Baum, Spin-Strukturen und Dirac-Operatoren \"{u}ber
pseudoriemannschen Mannigfaltigkeiten, Teubner-Verlag, Leipzig
(1981).
\bibitem{2}
H. Baum, Th. Friedrich, R. Grunewald, I. Kath,  Twistors and
Killing spinors on Riemannian manifolds, Teubner,
Leipzig/Stuttgart (1991).
\bibitem{3}
D. Bleecker,  Gauge theory and variational principles,
Addison-Wesley, Mass. (1981).
\bibitem{4}
J.P. Bourguignon, P. Gauduchon, Spineurs, Op\'{e}rateurs de Dirac
et Variations de M\'{e} triques,  Commum. Math. Phys. 144 (1992)
581-599.
\bibitem{5} Th. Friedrich,  Solutions of the Einstein-Dirac equation on Riemannian 3-manifolds
with constant scalar curvature, J. Geom. Phys. 36 (2000) 199-210.
\bibitem{6}
Th. Friedrich, E.C. Kim,  The Einstein-Dirac equation on
Riemannian spin manifolds, J. Geom. Phys. 33 (2000) 128-172.
\bibitem{7}
E. C. Kim,  A local existence theorem for the Einstein-Dirac
equation, J. Geom. Phys. 44 (2002) 376-405.
\bibitem{8}
M. Wang, Parallel spinors and parallel forms, Ann. Glob. Anal.
Geom. 7 (1989) 59-68.
\end{Literature}

\end{document}